\documentclass[12pt]{amsart}
\usepackage{graphicx}
\usepackage{lipsum}
\usepackage{xcolor}
\usepackage{booktabs}

\usepackage{amsmath,amsthm,amsfonts,amssymb}
\usepackage{amsmath, amsthm, amssymb, amscd, amsfonts}
\usepackage[all]{xy}
\usepackage{amssymb, amsthm, amsmath}
\usepackage[english]{babel}
\usepackage{amsfonts}

\usepackage{color}
 
\DeclareMathOperator{\id}{id}
\DeclareMathOperator{\Hom}{Hom}

\newcommand{\ra}{\rightarrow}

\newcommand{\C}{\mathbb C}

\newcommand{\ot}{\otimes}
\newcommand{\mtc}{\mathcal}
\newcommand{\lam}{\lambda}
\newcommand{\lb}{\label}
\newcommand{\Lam}{\Lambda}

\newcommand{\al}{\alpha}
\newcommand{\eps}{\epsilon}
\newcommand{\sub}{\subsection}

\newcommand{\D}{\Delta}

\newcommand{\lh}{\leftharpoonup}
\numberwithin{equation}{section}
\newtheorem{lm}[equation]{Lemma}

\theoremstyle{plain}
\newcommand{\un}{\mathbf{1}}
\newtheorem{thm}[equation]{Theorem}
\newtheorem{prop}[equation]{Proposition}
\newtheorem{defn}[equation]{Definition}
\newtheorem{cor}[equation]{Corollary}
\newtheorem{rem}[equation]{Remark}
\newcommand{\ch}{\chi}
\newcommand{\mtr}{\mathrm}

\numberwithin{equation}{section}
\newcommand{\ncm}{\newcommand}
\ncm{\np}{\newpage}
\ncm{\ebl}{\end{thebibliography}}
\ncm{\bbl}{\begin{thebibliography}}
\ncm{\chd}{_{ _{\ch}}}
\ncm{\ald}{_{ _{\al}}}
\ncm{\cP}{\mathcal{P}}
\ncm{\ei}{e_i}
\ncm{\eij}{e_{i,\;j}}
\ncm{\bt}{\begin{thm}}
\ncm{\bdef}{\begin{defn}}
\ncm{\edf}{\end{defn}}
\ncm{\et}{\end{thm}}
\ncm{\bc}{\begin{cor}}
\ncm{\bl}{\begin{lem}}
\ncm{\el}{\end{lem}}
\ncm{\bpf}{\begin{proof}}
\ncm{\epf}{\end{proof}}
\ncm{\ec}{\end{cor}}
\ncm{\ord}{\mtr{ord}}
\ncm{\er}{\end{rem}}
\ncm{\br}{\begin{rem}}
\ncm{\bn}{\begin}
\ncm{\bp}{\begin{prop}}
\ncm{\ep}{\end{prop}}
\newcommand{\tet}{\theta}

\newcommand{\rdual}[1]{{#1}^*}
\ncm{\bd}{\begin{document}}
\ncm{\ed}{\end{document}}
\ncm{\beq}{\begin{equation}}
\ncm{\beqn}{\begin{equation*}}
\ncm{\eeq}{\end{equation}}
\ncm{\eeqn}{\end{equation*}}
\ncm{\bea}{\begin{eqnarray}}
\ncm{\eea}{\end{eqnarray}}
\ncm{\beanon}{\begin{eqnarray*}}
\ncm{\eeanon}{\end{eqnarray*}}
\ncm{\ek}{\eps|_K}\ncm{\diez}{\#}
\ncm{\co}{\mtc{O}}
\ncm{\bwt}{\bowtie}
\newcommand{\xra}{\xrightarrow }
\newcommand{\ovr}{\overline}
\ncm{\cC}{\mtc{C}}\ncm{\cc}{\mtc{C}}
\ncm{\cX}{\mtc{X}}
\ncm{\wt}{\widetilde}
\ncm{\sg}{\sigma}\ncm{\Rep}{\mathrm{Rep}}
\ncm{\Irr}{\mathrm{Irr}}\ncm{\X}{\mathcal{X}}
\ncm{\cA}{\mathcal{A}}
\ncm{\HKer}{\mtr{HKer}}
\ncm{\LKER}{\mtr{LKer}}
\ncm{\aad}{\mtr{ad}}
\ncm{\Dr}{\mtr{D}}
\ncm{\cD}{\mathcal{D}}\ncm{\cd}{\mathcal{D}}\ncm{\ce}{\mathcal{E}}
\ncm{\irr}{\mtr{Irr}}
\ncm{\G}{\mathcal{G}}
\ncm{\Dc}{\mtc{D}}
\ncm{\E}{\mtc{E}}
\ncm{\fp}{\mtr{FPdim}}
\ncm{\Vc}{\mtr{Vec}}
\ncm{\cK}{\mtc{K}}
\ncm{\cM}{\mtc{M}}
\ncm{\cE}{\mtc{E}}
\ncm{\cS}{\mtc{S}}
\ncm{\End}{\mtr{End}}
\ncm{\cop}{\mtr{cop}}
\ncm{\op}{\mtr{op}}
\ncm{\chr}{character }\ncm{\ck}{\mtc{K}}
\ncm{\bw}{\bwt}
\ncm{\cY}{\mtc{Y}}
\ncm{\hker}{\mtr{HKer}}
\ncm{\bx}{\boxtimes}
\ncm{\blue}{\textcolor[rgb]{.00, .00, 1.00}}
\ncm{\red}{\textcolor[rgb]{1.00, .00, .00}}
\ncm{\green}{\textcolor[rgb]{.00, 1.00, .00}}
\ncm{\bne}{\begin{enumerate}}
\ncm{\ene}{\end{enumerate}}
\ncm{\lker}{\mtr{LKer}}
\ncm{\md}{\medbreak}
\ncm{\rep}{\Rep}\ncm{\ind}{\mtr{ind}}
\ncm{\mdn}{\md\noindent}
\newcommand{\kk}{\Bbbk}
\bd
\author{Sebastian Burciu}
\address{Inst.\ of Math.\ ``Simion Stoilow" of the Romanian Academy, Research Team 5, P.O. Box 1-764, RO-014700, Bucharest, Romania}
\thanks{This work was supported by a grant of the Romanian National Authority for Scientific Research, CNCS UEFISCDI, project number PN-II-RU-TE-2012-3-0168}
\title{ On an analogue of a Brauer Theorem for fusion categories}
\maketitle
\date{\today}
\bn{abstract}In this paper we prove an analogue of Brauer's theorem for faithful objects in fusion categories. Other notions, such as the order and the index associated to  faithful objects of fusion categories are also discussed. We show that the index of a faithful simple object  of a fusion categories coincides to the order of the universal grading group of the fusion category.
\end{abstract}

\section{Introduction}

In representation theory of finite groups a $\mathbb C G$-module $M$ is called {\it faithful} if its kernel $\ker_{\mathbb C G}(M)$ is trivial. A celebrated theorem of Brauer states that in this case any other representation of $\mathbb CG$ can be found as a constituent of at least one tensor power $M$. The goal of this paper is to generalize this result to fusion categories with commutative Grothendieck rings. We should also mention the fact that in the literature this theorem is also called Burnside-Brauer theorem.

Grothendieck rings of tensor and fusion categories were recently intensively studied by various authors. 
In this paper we give some new properties for the Grothendieck rings of  the fusion categories  generated  by a single  object. We say that an object $X \in \co(\cc)$ is {\it faithful} if the fusion subcategory $\cc(X)$ generated by $X$ coincides with the whole category $\cc$. Faithful objects were considered previously in the literature, see for example \cite{nfaith}. One of the main goals of this paper is to formulate a result analogue to Brauer' s theorem for faithful representations.

Recently, Brauer's theorem for representations of finite groups was generalized to Hopf algebras by several authors (\cite{KSZ, PQ,  Bker}). The starting point of all these generalizations was the results obtained by Rieffel in \cite{Riburn}. In all these references the notion of a faithful module is considered without defining the notion of kernel of a representation. More recently, the author defined the notion of left kernel and right kernel for modules over Hopf algebras.

In this paper we propose a notion of {\it kernel} for objects in arbitrary fusion categories.  We show that an analogue of Brauer's theorem holds in the case of fusion categories with a commutative Grothendieck ring.

In \cite{nfaith} the author has shown that if a fusion category has a faithful object then the universal grading group of $\cc$ is cyclic. In the present paper we show that the order of this cyclic group equals the {\it index} of the object $X$.  Analogue to \cite{KSZ}, one can define the  index of any object of $\cc$ as the index of the imprimitivity of the matrix associated to the operator $l_{[X]}$ of left multiplication by $[X]$ on the Grothendieck ring $\ck_0(\cc)$. For more details see Section  \ref{index}. 


This paper is organized as follows. In Section \ref{prelim} we recall the basics  from the theory of fusion categories that are needed through the paper. In the next section we introduce the notion of kernel of an object in a fusion category and prove an analogue of Brauer's theorem for fusion categories. We also introduce the notion of the center of a character and present some properties of the index of a character. In Section \ref{index} we discuss the notion of index of an object and also introduce the related concept of center of an object of a fusion category.  As an example, in Section \ref{four} we revisit the case of semisimple Hopf algebras and we show how the new notion of kernel relates to the previous notion of left (right) kernel introduced by the author in \cite{gmj}. Applications of Brauer's theorem to the case of a modular fusion category $\cc$ are presented in the last section of the paper.

All fusion categories and algebras in this paper are considered over an algebraically closed base field $\kk$ of characteristic zero.

\section{Preliminaries}\lb{prelim}
As usually, by a fusion category we mean a {\it $\kk$-linear semisimple rigid tensor category} $\cc$ with finitely many isomorphism classes of simple objects, finite
dimensional spaces of morphisms, and such that each simple object $S$ is scalar (that is, $\End_{\cc}(S) = \kk$) and the unit object $\un$ of $\cc$ is simple. We
refer the reader to \cite{eno} for basics on fusion categories. 
We denote by $\co(\cc)$ the class of all objects of $\cc$ and by $\irr(\cc)$ the set of isomorphism classes of simple objects of $\cc$.

A {\it fusion subcategory} of a given category is a full monoidal replete subcategory which is also fusion category. In a fusion category $\cc$, the left and right duals $X^{*}$ and $\rdual{X}$ of an object $X$ are isomorphic (but it is still not known whether there always exists a sovereign structure, that is a natural monoidal isomorphism between the two duals).

Recall that the Grothendieck group $K_0(\mtc{C})$ of a fusion category $\cc$ is the free $\mathbb{Z}$-module with basis given by the isomorphism classes of simple objects of $\cc$. 
It is a well known result  now (see also \cite{eno}) that $\ck_0(\cc)$ is a semisimple {\it based ring} with the ring structure endowed from the monoidal structure of $\cc$. 
By abuse of notation, if no confusion,  we sometimes may  write $X \in K_{0}(\cc)$ instead of $[X]\in \ck_{0}(\cc)$.
\subsection{Frobenius-Perron dimension of objects} If $A$ is a matrix with nonnegative entries then, by Frobenius-Perron theorem, the matrix $A$
has a positive eigenvalue $\lambda$ which has the biggest absolute
value among all the other eigenvalues of $A$ (see \cite{F}). The eigenspace
corresponding to $\lambda$ has a unique vector with all entries
positive. The eigenvalue $\lambda$ is called the {\it principal value} of $A$ and will also be denoted by $\fp(A)$. The corresponding positive eigenvector of $\lam$ is called {\it the principal vector} of $A$. Also the eigenspace of $A$ corresponding to $\lam$ is called the {\it principal eigenspace} of  the matrix $A$.

Let $\cc$ be a fusion category and $\ck(\cc):=\kk\ot_{\mathbb{Z}}\ck_0(\cc)$ be its Grothendieck ring. For an object $X\in \co(\cc)$ let $l_{[X]}$  be the operator of left multiplication by $[X]$ on $\ck(\cc)$. Also let $a(X)$ be the matrix associated to the operator $l_{[X]}$ with respect to the basis of $\ck(\cc)$ given by the simple objects of $\cc$. Then by definition one has that $\fp(X):=\fp(a(X))$ (see also \cite{eno}).  We denote by $\fp():\ck(\cc)\ra \kk$ the unique ring homomorphism for which $\fp(X)>0$ for any simple object $X$ of  $\cc$, see \cite{eno}. 
We also denote by $X(\cc)$ the set of all unitary ring homomorphisms $\fp():\ck_{0}(\cc)\ra \kk$.
\subsection{The symmetric associative bilinear form on $\ck_0(\cc)$.} Let $\cc$ be a fusion category and as above $\ck(\cc):=\mathbb{C}\ot_{\mathbb{Z}}\ck_0(\cc)$ be its Grothendieck ring.

It is well known that $\ck(\cc)$ is a semisimple $\kk$-algebra. Since $\kk$ is algebraically closed we may suppose that 
\beqn
\ck(\cc)\simeq M_{n_{0}}(\kk)\times M_{n_{1}}(\kk)\cdots \times M_{n_{s}}(\kk)
\eeqn
Let $E_0, E_1, \cdots, E_{s}$ be a complete set of primitive central idempotents of $\ck(\cc)$ and let $\mu_0, \cdots, \mu_s$ be the corresponding irreducible characters of $\cK_0(\cc)$. Without loss of generality we may assume that $E_0=\frac{1}{\fp(\cc)}\mtc{R}_{\cc}$ where $\mtc{R}_{\cc}$ is the virtual regular element of $\cc$ defined by \beq\mtc{R}_{\cc}=\sum_{X\in \Lam_{C}}\fp(X)[X].\eeq
Then the corresponding character of $E_{0}$ is  $\mu_0=\fp()$ given by $[X]\mapsto \fp(X)$ and $n_{0}=1$.
\md  Denote by $m_\cc$  the $\kk$-bilinear form on the Grothendieck ring $\ck(\mtc{C})$ defined on the  generators $X, Y\in \irr(\cc)$ by $m_{\cc}(X,\;Y) = \dim_\kk \Hom(X,Y)$ and then extended linerly on both arguments. 

The bilinear form
$m_\cc$ has the following properties:
\begin{enumerate}
\item symmetry: $m_\cc(x,y) = m_\cc(y,x)$;
\item adjunction property: $m_\cc(x,y\,z) = m_\cc(y^*\,x,z) = m_\cc(x\, z^*,y)$
\end{enumerate}
Define also $B:\ck_0(\cc)\ot \ck_0(\cc)\ra \kk$ by $B(x, y):=m_{\cc}(1, xy)$. Then it is easy to check that $B$ is a symmetric associative nondegenerate bilinear form on $\ck(\cc)$. It follows that there are scalars $b_{i}\in \kk^{*}$ such that 
\beq\label{m2'}
B(x, y)=\sum_{i=0}^sb_i\mu_i(xy)
\eeq
Moreover  the dual bases equation imply that
\beq\lb{db}
\sum_{i=1}^s\frac{1}{b_i}E^i_{jk}\ot E_{kj}^i=\sum_{X\in \Lam_{\cc}}[X]\ot [X^*]
\eeq
where $E^{i}_{jk}$ are the elementary matrix entries in $M_{n_{i}}(\kk)$. 
Note that $b_{i}\neq 0$ since the form is nondegenerate.

It follows from \cite{ostrik4} that $\frac{1}{b_i}=f_{\mu_i}$ is the formal codegree of the irreducible representation of $\ck(\cc)$. Therefore $f_{\mu_i}$ is a real positive number.
\md
Note that the above equalities also show that:
\beq\lb{m1}
m_{\cc}([X], [Y])=B([X],[Y]^*)=\sum_{i=1}^sb_i\mu_i([X][Y]^*)
\eeq
for any $X, Y\in \co(\cc)$.
\section{Brauer's theorem and faithful objects in fusion categories}\label{brfaith}
\subsection{Frobenius-Perron eigenvalues}
Let $\cc$ be a fusion category with the Grothendieck ring $\ck_0(\cc)$ commutative. Let $\psi:\ck(\cc)\ra \kk$ be a linear character and $E\in \ck(\cc)$ be its primitive central idempotent. Then since $E$ is an eigenvector for $l_{[X]}$ it follows that $|\psi([X])|\leq \fp(X)$ for any object $X \in \co(\cc)$. Indeed $[X]E=\psi([X])E$ and therefore $E$ is an eigenvector for the operator $l_{[X]}$, 
left multiplication by $[X]$.

\bn{defn}
Define $\ker_{\cc}(\psi)$ as the set of all simple objects $X$ such that $\psi([X])=\fp(X)$.
\end{defn}
\br
It is known that $\psi(X) \in \mathbf{Q}(\xi)$ for some root of unity $\xi$, see  \cite[Corollary 8.53]{eno}. 
\er

\bp \lb{preigenv}Assume  that the Grothendieck ring $\ck_0(\cc) $ is commutative.
If $X$ and $Y$ have a principal eigenvector $V\in \ck_0(\cc)$ then any simple object of the fusion subcategory generated by $X $ and $Y$ have $V$ as a principal eigenvector.
\ep
\bpf
Suppose that $V=\sum_{i=0}^s \mu_i(V)E_i$. Since $XV=\fp(X)V$ it follows that $
\fp(X)=\mu_i(X)$ for all $i$ with $\mu_i(V)\neq 0$. Therefore it is enough to consider the case $V=E_i$ some idempotent of $\ck(\cc)$.
\md Suppose that $X, Y \in \ker_{\cc}(\psi)$ and let 
\beq
X\ot Y=\oplus_{Z}N^Z_{X,Y}Z
\eeq 
If $|\psi([X])|= \fp(X)$ and  $|\psi([Y])|=\fp(Y)$ then
\beqn
|\psi(X\ot Y)|=|\oplus_{Z}N^Z_{X,Y}\psi(Z)|\leq \oplus_{Z}N^Z_{X,Y}\fp(Z)=\fp(X\ot Y)
\eeqn
Since $\psi$ is an algebra homomorphism it follows that $\psi(X\ot Y)=\fp(X)\fp(Y)$ and the previous inequality shows that $\psi(Z)=\fp(Z)$ for all $Z$ with $N^Z_{XY}\neq 0$.
\epf

Similalry one can prove the folllowing:
\bc\lb{const}Assume that the Grothendieck ring $\ck_0(\cc) $ is commutative.
If $[X]V=\xi\fp([X])V$ for some object $X\in \co(\cc)$ and some root of unity $\xi$ then $[Y]V=\xi\fp([Y])V$ for any simple subobject $Y$ of $X$.
\ec
Proposition \ref{preigenv} implies:
\begin{lm}\lb{fc}
The set $\ker_{\cc}(\psi)$ is closed under tensor products and therefore it spans a fusion subcategory of $\cc$.
\end{lm}

\bn{exmp}\lb{grp}
Let $G$ be a finite group and $\cc=\rep(G)$. Then $\ck(\cc)=C(G)$, the character ring of $G$. Any ring homomorphism of $C(G)$ corresponds to a conjugacy class $\mtc{X}$ of $G$. It is given by  $\widehat{\mtc{X}}$ which is the characteristic function of $\mtc{X}$ on $G$. Then clearly $\ker_{\cc}(\widehat{\mtc{X}})=\ker_{kG^*}(g)$ where $g\in \mtc{X}$ is an arbitrary group element of $\mtc X$.
\end{exmp}

\bn{defn}\label{kernc}
For any object $X\in \co(\cc)$ define 
\beq
\ker_{\cc}([X]):=\{\psi\in \irr(\ck(\cc))\;|\;[X]\in \ker_{\cc}(\psi)\}
\eeq
\end{defn}
Note that for any object $X\in \co(\cc)$ one has that $\fp()\in \ker_{\cc}(X)$.\mdn
We say that an object $X$ has {\it a trivial kernel} if $\ker_{\cc}([X])=\{\fp()\}$. Next theorem can be regarded as a generalization of Brauer's Theorem.
\bt\lb{brauer}
Let $\cc$ be a fusion category with commutative Grothendieck ring.
If $X$ is an object of $\cc$ with trivial kernel then any other object of $\cc$ is a subobject  of some tensor power of $X$.
\et
\bpf Using Equation \eqref{m2'} the proof follows the lines of the proof from the classical case of group representations, see for example \cite [Theorem 19.10]{liebeck}. Note that if $\ck(\cc)$ is a commutative ring then Equation \eqref{m2'} implies that
\beq\lb{m2}
B(x, y)=\sum_{i=0}^sb_i\mu_i(x)\mu_i(y)
\eeq
for any $x,y \in \mtc{K}(\cc)$.

Suppose that \beqn
[X]=\sum_{i=0}^s\mu_i([X])E_i
\eeqn
and let $Y \in \irr(\cc)$ with
\beqn
[Y]=\sum_{i=0}^s\mu_i([Y])E_i
\eeqn
Since $\ker_{\cc}(X)=\{\fp()\}$ it follows that $\mu_i([X])\neq \fp(X)$ for $i \geq 1$.
If $Y$ is not a direct summand in any tensor power of $X$ then by Equation \eqref{m2} one has that:
\beq\lb{multip}
\sum_{i=0}^s\mu_i([X])^n\mu_i([Y])=0.
\eeq
for any $n \geq 1$. Denote by $\mtc A_i$ the set  of all indices $ 0\leq j \leq r$ with $\mu_j([X])=\mu_i([X])$. It follows that $\mtc A_0=\{0\}$.
\mdn
Then by taking a van der Monde determinant, Equation \eqref{multip} implies that each  $\sum_{j\in \mtc{A}_i}{n_j\mu_j(Y)}=0$. Note that this is impossible for $i=0$ since $n_0\fp(Y)\neq 0$. 
\epf
\ncm{\ca}{\mtc{A}}\ncm{\cz}{\mtc{Z}}

\subsection{On the Grothendieck ring of $\cc(X)$}
Let $\cc$ be a fusion category with commutative Grothendieck ring $\ck(\cc)$. We define $\cd:=\cc(X)$. Clearly $\ck(\cd)\subset \ck(\cc)$ and therefore one can write the primitive central idempotent $\frac{r_{\cd}}{\fp(\cd)}$ as sum of central primitive idempotents of  $K(\cc)$:
\beqn
\frac{r_{\cd}}{\fp(\cd)}=\sum_{j \in \mtc{A}_X}E_j
\eeqn
for a subset $\ca_X \subseteq \{0, 1, \cdots , s\}$.
\bp\lb{retk}
With the above notations one has that 
\beq
\ker_{\cc}(X)=\{\mu_j\;|\;j \in \ca_X\}
\eeq
\ep\newcommand{\dd}{$\;}

\bpf
Let $F_0, F_1,\cdots, F_r$ be the primitive idempotents of $\ck(\cd)$ and $\ch_0, \cdots, \ch_r$ be their associated irreducible characters. Moreover as above one may suppose that $F_0=\frac{r_{\cd}}{\fp(\cd)}$ and consequently $\ch_{0}=\fp()$.

We look at the restrictions of the characters $\mu_{i}$ at $K(\cd)$. Then there is a surjective function $f:\{0, \cdots ,s\}\ra \{0, \cdots, r\}$ such that ${\mu_j}{|_{\ck(\cd)}}=\ch_{f(j)}$ for all $0 \leq j \leq s$.

With the above notations it will be shown that $f^{-1}(0)=\ca_{X}$.
Indeed, since $X$ is a faithful object of $\cd$ it follows that $\ker_{\cd}(X)=\{\mu_0\}$.
On the other hand note that $\ch_j([X])={\ch_j}{|_{\ck(\cd)}}([X])=\mu_{f(j)}([X])$. Thus $\ch_j \in \ker_{\cc}(X)$ if and only if $\mu_{f(j)}\in \ker_{\cd}(X)$, i.e. $f(j)=0$.
\epf
One can also formulate the converse of the  Brauer's theorem in categorical settings:
\bp
If a fusion category $\cc$  with the Groethendieck ring $\ck(\cc)$ commutative is generated by a single object $X$ then $\ker_{\cc}(X)=\{\fp()\}$ 
\ep
\bpf
Suppose that $\mu_i([X])=\fp(X)$ for some $1 \leq i\leq s$. Then Lemma \ref{fc} implies that $\mu_i(Y)=\fp(Y)$ for any object $Y \in \co(\cc)$ and therefore $\mu_i=\fp()$.
\epf
\bn{exmp}
Let $G$ be a finite group and $\cc=\rep(G)$. Then as in Example \ref{grp} one has $\ck(\cc)=C(G)$. It follows that $\ch \in \Irr(G)$ is faithful in the above sense if and only if the kernel of $\ch$ is trivial, i.e. it is faithful in the classical sense.
\end{exmp}

\section{The index of a simple object and its center}\label{index}

The second part of Frobenius-Perron theorem states that for a nonnegative indecomposable matrix $A$ any other eigenvalue $\mu$ of $A$ such that $|\mu|=\fp(A)$ has a unique (up to a scalar) corresponding eigenvector. 
The number of such eigenvalues $\mu$ with $|\mu|=\fp(A)$ is called the {\it imprimitivity index} of $A$. Moreover, if $\xi$ is a primitve root of unity of order equal to the index $A$ then $\mu=\xi^{i}\fp(A)$ for some $0\leq i \leq \ind(A)-1$.
\md
Let $\cc$ be a fusion category.
One can define the {\it index} of an object as the index of the imprimitivity of the matrix $a(X)$ associated to the operator $l_{[X]}$ on $\ck(\cc)$ with respect to the canonical basis of $\ck(\cc)$ given by the simple objects of $\cc$. 
\md
Similarly to \cite{KSZ} one can show that an object $X \in \co(\cc)$ is faithful if and only if the associated matrix $a(X)$ of $l_{[X]}$ is indecomposable. Recall \cite{F} that a $n \times n $-matrix $A = (a_{ij})$ is called {\it decomposable} if it is
possible to find a decomposition of $\{1, 2, \cdots, n\}= M \cup N $ of into disjoint
nonempty sets $M$ and $N$ such that $a_{ij} = 0 $ whenever $i \in M$ and $j \in N$ . Otherwise the matrix  $A$ is called {\it indecomposable. }

\bt\label{indexdiv}  Let $\cc$ be a fusion category. 
 Suppose that a simple object $Y$ of $\cc$ is a constituent of some tensor powers $X^{\ot m}$ and $X^{\ot n}$ of $X$ with $m, n \geq 0$. Then $m-n$ is divisible by $\ind(X)$.
\et
\bpf
One may assume that $\cc=\cc(X)$. Let $\xi$ be a primitive root of unity of order $\mtr{index}(X)$. Let  also $v \in K(\cc)$ be the unique (up to a scalar) eigenvector for the linear operator $l_{[X]}$ corresponding to the eigenvalue $\xi \fp(X)$.  It follows that for any other object $Z\in \co(\cc)$ one has that $v[Z]$ is also an eigenvector for $l_{[X]}$ corresponding to the same eigenvalue $\xi \fp(X)$. Therefore, by the above discussion there is a scalar $\gamma([Z])\in \kk$ such that
$
v[Z]=\gamma([Z])v
$
for all $Z \in \co(\cc)$.
\md
Clearly $\gamma:\ck(\cc)\ra \kk$ is an algebra homomorphism. Since $\ck(\cc)$ is a semisimple algebra it follows that $v$ is a central element of $\ck(\cc)$. 
Therefore one can write that
\beqn
v=\sum_{i \in \ca}v_iE_i.
\eeqn
for some subset $\ca\subseteq \{0,1, \cdots, s\}$ and $v_{i}\neq 0$ for $i \in \ca$.
\md
Since $v$ is unique up to a scalar it follows that $[X]E_i\neq 0$ for any $i \in \ca$. Therefore $\ca$ should be a set with one element. We may suppose without loss of generality that $\ca=\{1\}$. Therefore one can write $[X]^mE_1=\xi^m\fp(X)^mE_1$.  Then Lemma \ref{const} implies that $\mu_1([Y])=\xi^m\fp(Y)$ since  $Y$ is a constituent of $X^{\ot \; m}$. On the other hand the same lemma implies that  $\mu_1([Y])=\xi^n\fp(Y)$ since $Y$ is also constituent of $X^{\ot \; n}$. Thus $m-n$ is divisible by the index of $X$.\epf

Note that in \cite[Theorem 4.1]{nfaith} it was proven that the universal grading group $U(\cc(X))$ is cyclic. Next theorem gives a precise description for the order of this cyclic group. It can also be seen as a generalization of the result \cite[Proposition 4.5]{KSZ}.
\bt\label{order}
The universal grading group of $\cc(X)$  is isomorphic to $\mathbb{Z}_n$ where $n=\ind(X)$.
\et
\bpf
As we already mentioned above we know that $U(\cc(X))$ is a cyclic group. By  \cite[Proposition 3.1]{ng} it is enough to show that one has $\mtr{Aut}_{\ot}(\id_{\cc(X)})\cong \mathbb{Z}_n$. Note that in order to give a $\phi \in \mtr{Aut}_{\ot}(\id_{\cc(X)})$ it is enough to give for each $X \in \irr(\cc)$ a scalar $s(X)$ such that $\phi_X=s(X)\id_X$. Then one can to define
\beq
\phi_Y=s(X)^m\id_Y
\eeq
for any object $Y \in \co(\cc(X))$ with $Y$ a constituent of $X^{\ot m}$.
By previous theorem it follows  that $\phi$ is well defined. Indeed  if $Y$ is a direct summand in $X^{\ot m}$ and $X^{\ot n}$ then $m-n$ is divisible by $\ind(X)$.  Moreover $s(X)$ should be a root of unity of order $\ind(X)$ since $X^{\ind(X)}$ contains $1$ and $\phi_{1}:1\ra 1$ is identity. 
Conversely, any root of unity of order $\ind(X)$ determines a unqiue element $\phi \in Aut_{\ot}(\id_{\cc(X)})$. \epf
\subsection*{Construction of the universal grading of $\cc(X)$}
Let $\xi$ be a primitive root of unity of order $\ind(X)$.  As in the proof of  Theorem \ref{indexdiv} let  also $E_1$ be the central idempotent corresponding to the eigenvalue $\xi\fp(X)$ where $\xi$ is a primitive root of unity of order $\mtr{ind}(X)$.
Thus
\beqn
E_1[Y]=\xi^m\fp(Y)E_1
\eeqn
for all objects $Y$ that are constituents of $X^{\ot m}$.  \md
For all $0 \leq i \leq n-1$  define the full abelian subcategory $\cd_i$ of $\cc(X)$ by
\beqn
\cd_i:=\{Y\in \irr(\cc(X))\;|\mu_1([Y])=\xi^i\fp([Y])\}.
\eeqn
Then it can easily check that $\cd_i \ot \cd_j$ is mapped to $\cd_{i+j}$.
This shows that the decomposition $\cc=\oplus_{i=0}^{n-1}\cd_i$
is a grading of $\cc$. By the above arguments this grading should coincide with the universal grading of $\cc(X)$.
\subsection{The order of an object}
As explained in \cite{KSZ, nfaith} for any object $X \in \co(\cc)$ there is a smallest integer $n \geq 0$ such that $m_{\cc}(1_{\cc}, X^{\ot \; n})>0$, i.e $X^{\ot \; n}$ contains the unit object $1_{\cc}$. Then $n$ is called {\it the order of $X$} and denoted by $o(X)$.
This also shows that $X^*$ is also subobject of some tensor power of $X$.  If $X$ is a self dual object clearly its order is $2$.
\bc
The index of $[X]$ divides the order of $X$.
\ec
\bpf
If $X^m$ contains the unit element then under the above grading it follows that $\xi^m=1$ and therefore $\ind(X)|\ord(X)$.
\epf
\subsection{The center of a an object}
\bn{defn} Let $\cc$ be a fusion category over $\kk$ with a commutative Groethendieck ring and $X$ be an object of $\cc$.
Define 
\beq
\cz_{\cc}(X)=\{\mu_i\in \Irr(\ck(\cc))\;|\;\;|\mu_i([X])|=\fp(X)\}
\eeq
\end{defn}
 Note that similarly to \cite[Proposition 3.1]{ws3} or \cite[Proposition 10]{NR'} one has that $\psi([X^{*}])=\overline{\psi([X])}$ for any simple object $X$ of $\cc$ and any character $\psi:\ck(\cc)\ra \kk$. This implies that for any simple object $\ker_{\cc}(X\ot X^{*})=\cz_{\cc}(X)$ and  therefore
\beqn
\ker_{\cc}(\mtr C_{ad})=\bigcap_{X \in \irr(\cc)}\cz_{\cc}(X)
\eeqn
whre $\mtr C_{ad}:=\bigoplus_{X\in \irr(\cc)}X\otimes X^{*}$.

For any $0 \leq i \leq \ind(X)-1$ let $\mu_i$ be the linear character of $\ck(\cc(X))$ corresponding to the unique (up to a scalar) eigenvector of $l_{[X]}$ with eigenvalue $\xi^i\fp(X)$ in $\ck(\cc(X))$.
\bc
With the above notations one has that:
\beq \cz_{\cc(X)}(X)=\{\mu_0, \cdots ,\mu_{\ind(x)}\}. \eeq
\ec
\section{Examples from semisimple Hopf algebras}\label{four}
 Let $A$ be a  Hopf algebra over a field $\kk$. We use the standard Hopf algebra notations that can be found for example in \cite{montg}. We denote by $\rep(A)$ the tensor category of finite dimensional representations of $A$. If $A$ is a semisimple Hopf algebra then it is known that  $\rep(A)$ is a fusion category.

Let $M$ be a finite dimensional left $A$-module. Then the left kernel $\lker_A(M)$ of a finite dimensional representation $M$ of $A$ is defined by
\beq
\lker_A(M)=\{a \in A\;|\; a_1\ot a_2m=a\ot m\;\text{for all}\; m \in M\}.
\eeq
It follows (see \cite{gmj}) that $\lker_A(M)$ is the largest left coideal subalgebra of $A$ that acts trivially on $M$.
Recall that a left coideal is a subspace $L$ of $A$ with the property that $\Delta(L)\subseteq A\otimes L$. A {\it left coideal subalgebra } is a left coideal which also is a subalgebra of $A$.

For a subspace $S$ of $A$ we denote by $S^{+}$ the subspace $S \cap \mtr{ker} \eps$.
If $L$ is a normal coideal subalgebra then it is known that $AL^{+}$ is a Hopf ideal and therefore $A/AL^{+}$ is a quotient of Hopf algebra. Via the canonical projection $\pi_{L}:A\ra A//L$ one can view $\rep(A//L)$ as a fusion subcategory of $\rep(A)$.

Brauer's theorem for Hopf algebras, see \cite[Theorem 4.2.1]{gmj} states  that $\rep(A//\lker_A(M))$ is precisely the fusion subcategory of $\rep(A)$ generated by the object $M$.

\bp\lb{incl}\bn{enumerate}
\item
Suppose that $\cd$ and $\ce$ are fusion subcategories of  a fusion category $\cc$. If $r_{\cd}r_{\ce}=\fp(\cd)r_{\ce}$ then $\cd\subseteq \ce$.
\item
Let $L$ and $K$ be two normal coideal subalgebras of a Hopf algebra $A$. If $AL^+\subseteq AK^+$ then $L\subseteq K$.
\end{enumerate}
\ep
\bpf
i) The right coset equivalence $r^{\cc}_ {\ce}$ on  the set of isomorphism classes of simple objects $\irr(\cc)$ (see \cite{bubr}) implies that any element of $\cd$ is equivalent to the unit $1_{\cc}$ and therefore it is contained in $\ce$. 

ii) If $AL^{+} \subseteq AK^{+}$ then one has a canonical Hopf projection $\pi_{K, L}:A/AL^{+}\ra  A/AK^+$. Moreover the the composition of the Hopf projections $A\xrightarrow{\pi_{L}} A/AL^+\xrightarrow{\pi_{K, L}} A/AK^+$ coincides to $\pi_{K}$. Then following \cite{tkq} one has that $L=A^{\mtr{co} \pi_{L}}\subseteq A^{\mtr{co} \pi_{K,L}\circ\pi_{L}}=A^{\mtr{co} \pi_{K}}=K$.
\epf

\subsection{Central idempotents of the character algebra $C(A)$}
Let $A$ be a semisimple Hopf algebra and $D(A)$ be its Drinfeld double. Then $A$ is a $D(A)$-module  (see \cite{zind}) via the scion
\beq \lb{zhu}(f \bwt a).b=(a_1bS(a_2))\lh S^{-1}f\eeq
for all $f \in A^*$ and $a, b\in A$.

For a decomposition $A=V \oplus W$ of $A$ as a direct sum of two $D(A)$-submodules define a linear functional $p_V\in A^*$ such that $p_V(v)=\eps(v)$ and $p_V(w)=0$ for all $v \in V$ and all $w\in W$.

Let now 
\beq
A = V_1 \oplus V_2 \cdots \oplus V_s
\eeq be a decomposition of $A$ as sum of simple $D(A)$-modules. Note that the simple $D(A)$-submodules of $A$ are the minimal normal left coideals of $A$. Recall that a left coideal  is called {\it normal}  if it closed under the left adjoint action i.e. $x_{1}MS(x_{2})\subseteq M$ for any $x \in M$.

By \cite[Theorem 15.3]{nrm} it follows that $\{p_{V_i}\}$ is a complete set of central primitive orthogonal idempotents of $C(A)$. Recall that the character algebra $C(A)$ is defined as ${K}_{0}(\rep(A))\ot_{\mathbb Z}\mathbb C$.
\mdn
\bt\label{relkern}
Let $A$ be a Hopf algebra with $C(A)$ a commutative ring (e.g. quasitriangular). Then  one has that:
\beq
\ker_{\rep(A)}(M)=\{\mu_V\;|\;V\subset \lker_A(M)\}
\eeq
\et
\bpf Note that if $C(A)$ is a commutative ring the \cite[Theorem 15.3]{nrm} implies that a linear basis of $C(A)$ is given by $p_{V}$ with $V\subset A$ a simple $D(A)$-submodule. Denote by $\mu_{V}:C(A)\ra k$ the corresponding algebra homomorphism. Also note that any normal left coideal subalgebra $L$ of $A$ can be regarded as a $D(A)$-submodule of $A$.

Suppose now that 
\beqn
\ch_M=\sum_{V\subset A}\mu_V(\ch_M)p_V
\eeqn
Then it follows that $\ch_M(v)=\eps(v)\mu_{V_{i}}(\ch_M)$ for any $v \in V_{i}$.

Suppose now that $V\subset \lker_A(M)$. Then $vm=\eps(v)m$ for all $m \in M$  since any element of the left kernel $\lker_{A}(M)$ acts trivially on $M$. Therefore $\ch_M(v)=\eps(v)\ch_M(1)$ for any $v \in V$. This implies that $\ch_Mp_V=\ch_M(1)p_V$ and therefore $\mu_V(\ch_M)=\ch_M(1)$. Hence $
\mu_V\in \ker_{\rep(A)}(\ch_M)$.

Recall that the regular character of a semisimple Hopf algebra $A$ coincides to the integral on the dual with the value at $1$ equal to the dimension of $A$. Using this,  by Lemma  \ref{retk} one has that
\beqn
\ker_{\rep(A)}(\ch_M)=\{ \mu_V \;|\;p_Vt=p_V\}
\eeqn
where $t \in (A//\lker_A(M))^*$ is the nonzero integral with $t(1)=1$. Let $\pi:A \ra A//L$ be the canonical Hopf projection.
Note that $p_Vt_{}=p_V$ if and only if $t(\pi(x))=\eps(x)$ for any $x \in V$. Indeed on the $D(A)$ submodule complement of $V$ one has that both $tp_V$ and $p_V$ vanishes. On the other hand on $x \in V$ one has that $p_Vt(x)=t(\pi(x))$ while $p_V(x)=\eps(x)$. \md Since $\pi(x_1)t(\pi(x_2))=t(\pi(x))\pi(1)$ one obtains that $t(\pi(x))=\eps(x)$ for any $x \in V$ if and only if $V^+\subset L^+$. The implies that 
\beq\label{krdsc}
\ker_{\rep(A)}(\ch_M)=\{\mu_V\;|\; V^+\subseteq AL^+\}
\eeq
Let $<V>$ denote the subalgebra of $A$ generated by $V$ inside $A$. Note that $<V>$ is also a left normal coideal subalgebra of $A$. 
In the situation of Equation \eqref{krdsc} one has $A<V>^+\subseteq AL^+$ and therefore $<V>\subseteq L$ by  Proposition \ref{incl}.
\epf
Following \cite{gmj} from the previous theorem one gets that:
\bc
Let $\cc=\Rep(A)$ be the category of finite dimensional  representations of a semisimple Hopf algebra $A$ and $M$ a finite dimensional representation of $A$. Then $M$ as an object of the fusion category $\Rep(A)$ is faithful if and only if $\lker_A(M)=k$.
\ec
\br
Remark that the proof of the previous theorem shows that 
\beq
\ker_{\rep(H)}(M)\supseteq \{\mu_V\;|V\subset \lker_{H}(M)\;\}
\eeq
even when the character ring $C(H)$ is noncommutative.
\er

\subsection{On the centre of a representation of a Hopf algebra} Let $A$ be a semisimple Hopf algebra and $\pi:A \ra \End_{k}(M)$ be a finite dimensional representation  of $A$. Recall that as in \cite{nfaith} one can define a notion of center by 
\beqn
LZ(\pi):=\lker_{A}(\pi \ot \pi^{*})
\eeqn
Then as in \cite[Proposition 3.5]{nfaith} $LZ(\pi)$ is the largest left coideal subalgebra of $A$ for which $\pi(a)=\lam(a)\pi(1)$ for a linear character $\lam:A\ra \kk$.

Note that similarly to the definition ofthe left kernel, see \cite{gmj}, one can also define $LZ(\pi)$ as the set of all elements $a \in A$ for which
$
\sum a_{1}\ot a_{2}m=\sum a_{1}\lam(a_{2}) \ot m
$
for some scalars $\lam(a_{2})\in \kk$. 

For a semisimple Hopf algebra we denote by $K(A)$ the largest central Hopf subalgebra of $A$. Moreover by \cite[Theorem 3.8]{ng} one has that $K(A)=k^{U(A)}$ where $U(A)$ is the universal grading group of $A$. We also denote by $A_{ad}$ the left adjoint $A$-module. Thus $A_{ad}=A$ as vector spaces and the module structure is given by $a.b=a_{1}bS(a_{2})$ for all $a,b \in A$. Since $A_{ad}$ generates the fusion subcategory $\rep(A//K(A)$ (see \cite[Proposition 18]{bd}) it follows that $\lker_{A}(A_{ad})=K(A)$. On the other hand it is not difficult to check by definition that $\lker_{A}(A_{ad})$ is the largest left coideal subalgebra contained in the center of $A$.
Thus, Brauer's theorem for Hopf algebras implies that $K(A)$ is also the largest central left (right) coideal subalgebra of $A$.
\bt Let $A$ be a semisimple Hopf algebra with a commutative character ring and $\ch$ be the character of a finite dimensional representation $M$  of $A$.
With the above notations one has the following
 \beqn LZ(\pi)/\lker(\pi)\subseteq K(A/\lker(\pi)).\eeqn Moreover $K(A/\lker(\pi))$ is the Hopf algebra of a cyclic group of order equal to the index of $\pi$. If $\pi$ is an irreducible character then one has equality above.
\et
\bpf
Let $L:=\lker_{A}(\pi)$. By Brauer's theorem, any irreducible representation of $A//L$ is a constituent of some tensor power. One has that $LZ(\pi)//L$ is a left coideal subalgebra contained in the center of   $A//L$ since any of its elements act as scalars on any irreducible representation of $A//L$.
\epf

Note that the above result generalizes \cite[Lemma 2.27]{Is} from group representations to Hopf algebras representations.

\section{The  case of a modular fusion category} In this section we apply our previous results to modular fusion categories, (see \cite{BaKi, Tu}). Recall that a {\it braided} tensor category $\cc$ is a tensor category equipped for all $X, Y \in \co(\cc)$ with natural isomorphisms $c_{X, Y}:X\ot Y \ra Y\ot X$ satisfying the hexagon axiom, see for example \cite{BaKi, js}.

A {\it twist} on a braided fusion category $\cc$ is a natural automorphism $\tet:\id_{\cc}\ra \id_{\cc}$ satisfying $\tet_{1}=\id_{1}$
and
$
\tet_{X\ot Y}=(\tet_{X}\ot \tet_{Y})c_{YX}c_{XY}.
$
A braided fusion category is called {\it premodular} or {\it ribbon} if it has a twist satisfying $\tet_{X^{*}} =\tet_{X}^{*}$ for all $X \in \co(\cc)$.

Recall that the entries of the $S$-matrix,  $S=\{s_{X,Y}\}$ of a premodular category are defined as the quantum trace $s_{X,Y}:=tr_{q}(c_{YX}c_{XY})$, see \cite{Tu}. Then it follows from \cite{proclond} (see also \cite[Lemma 6.5]{ng}) that $|s_{X,\;Y}|\leq \fp(X)\fp(Y)$ and  $s_{X,\;Y}=\fp(X)\fp(Y)$ if and only if \beq\lb{centr} c_{X,\; Y}c_{Y,\;X}=\mtr{id}_{XY}.\eeq
In the situation of Equation \eqref{centr} we say that $X$ and $Y$ {\it centralize each other}.  Morever {\it the centralizer } $\cd'$ of a fusion subcategory $\cd$ is defined as the full fusion subcategory of $\cd$ generated by all the objects of $\cc$ that centralize any object of $\cd$.  We say that two objects $X$ and $Y$ {\it projectively centralize} each other if $c_{X, Y}c_{Y, X}=\omega \fp(X)\fp(Y)$ for a root of unity $\omega\in \kk$.

A premodular category $\cc$ is called {\it  modular} if the above $S$-matrix is nondegenerate. Recall also that a fusion subcategory $\cd$ of a tensor modular category $\cc$ is called {\it nondegenerate } if $\cd\cap \cd'=\mtr{Vec}$. 
\md

Now let $\cc$ be a modular fusion category. For any $X \in \irr(\cc)$ the assignment $s_{X}([Y])=\frac{s_{X, Y}}{s_{X, 1}}$ extends linearly to a ring homomorphism $\ck(C)\ra \kk$, see \cite[Theorem 3.1.1]{BaKi}.  Conversely, any such ring homomorphism  is of the type $s_{X}$ for a simple object $X$  of $\cc$. Thus in the modular case one has 
\beq\label{moddesc}
X(\cc)=\{s_{X}\;|\;X \in \irr(\cc)\}.
\eeq

In this situation it follows that 
\beqn
\ker_{\cc}(s_{X})=\{Y\in \irr(\cc)\;|\; Y \;\text{centralizes}\; X\}=\cc(X)'.
\eeqn
and
\beqn
\cz_{\cc}(s_{X})=\{Y\in \irr(\cc)\;|\; Y \;\text{projectively centralizes}\;  X\}
\eeqn
By \cite[Proposition 6.7]{ng} it follows that 
\beqn
\cz_{\cc}(s_{X})=(\cc(X)')^{\mtr{co}}=(\cc(X)_{ad})'.
\eeqn
Therefore it follows from Theorem \ref{brauer} that in a modular category a simple object $X$ of $\cc$ is {\it faithful} if and only if $X$ does not centralize any other object than the trivial object. 

\md
For $S\subseteq X(\cc)$  we define $\widetilde{S}$ as the subset of simple objects of $\cc$ such that $s_{X}\in S$.
\bc
Suppose that $\cc$ is a single generated braided fusion category. Then $\widetilde {Z_{\cc}(X)}:=\langle Y\;|\; s_{Y}\in Z_{\cc}(X)\rangle$ coincides to group of invertible objects of $\cc$.
\ec

\br
It follows from \cite[Lemma 6.1]{ng}  that the simple objects $X_i$ that projectively centralize any other object of $\cc$ are precisely the invertible objects of $\cc$. Thus one can write that:
\beq
\bigcap_{i \in I}\widetilde{\cz_{\cc}(X_i)}=X(\cc).
\eeq
\er

\subsection{Relation with quantum doubles and kernel of the adjoint object}
Let $\cc$ be a spherical fusion category. Then it is  known that $\cz(\cc)$ is a modular tensor category see \cite[
Theorem 1.2]{mu3} or \cite[Theorem 2.15]{eno} with $\fp(\cz(\cc) ) = \fp(\cc)^{2}.$ Moreover   the forgetful functor $F:\cz(\cc) \ra \cc$  induces a surjective ring homomorphism $F_{*}:\ck(\cz(\cc)) \ra \ck(\cc)$, see \cite{eno}.

Now let $ E \in \irr(\ck(\cc))$. Then $ \ck(\cz(\cc))$ acts irreducibly on E via $F_{*}$. Let $\hat E \in \irr(\ck(\cz(\cc))) $ be the corresponding 1-dimensional representation of $ \ck(\cz(\cc))$. Thus by Equation \eqref{moddesc} there exists a unique $A_{E} \in \irr(\cz(\cc))$ such that $ \hat E = s_{A_{E}}$  and in addition
\beqn
f_{\hat E} =\frac{(\dim \cc)^{2}}{(\dim A_{E})^{2}}
\eeqn 
Moreover by \cite[Theorem 2.13]{o3} we have
\beq\label{des}
\dim A_{E}=\frac{\dim \cc}{f_{E}}
\eeq
and
\beqn
[I(1):A_{E}]=\dim E
\eeqn
Write $C_{ad}:=\bigoplus_{X\in \irr(\cc)}X\ot X^{*}$
\bp Let $\cc$ be a spherical fusion category. With the above notations it follows that
\beq
\ker_{\cc}(C_{ad})=\{F_{*}(s_{A})\;|\; A\in \mtr{Inv}(\cz(\cc)), [I(1):A]>0\}
\eeq
where $\mtr{Inv}(\cz(\cc))$ is the set of isomorphism classes of invertible objects of $\cz(\cc)$.
\ep
\bpf
Using Equation \eqref{db} it follows that $C_{ad}=\sum_{i=0}^{r}f_{i}E_{i}$. Moreover there is a simple object $A\in \irr(\cz(\cc))$ such that $\mu_{i}=F_{*}(s_{A})$. Moreover by Equation \eqref{des} one has that $f_{i}=\frac{\dim \cc}{\dim A}$. Thus $f_{i}=\fp(\cc)$ if and only if $\dim A=1$, i.e $A$ is an invertible object of $\cz(\cc)$.
\epf

\subsection*{Acknowledgements} The author thanks Dmitri Nikshych for fruitful conversations on paper's topic during the meeting ``Quantum Groups'' from Clermont-Ferrand September 2010. Part of this research was done  during a stay at the Erwin Schr\"odinger Institute, Vienna, in the frame of the Programme ``Modern Trends in Topological Quantum Field Theory'' in February 2014. The author also thanks the ESI and the organizers of the programme for their support and their very kind hospitality.
\bibliographystyle{amsplain}
\bibliography{hbb}
\ed\newpage
\bc Let $\cc$ be a modular tensor category.

1)If $\cc$ is single generated then $\cc$ is prime.

2) The number of generators of a modular tensor category is greater or equal then the number of components.
\ec
\bpf If
\beqn
\cc=\cd_{1}\boxtimes \cd_{2}
\eeqn
then $X$ is in exactly one of the components.
\epf
\bn{exmp}\label{singledg}
If $D(G)$ admits a faithful object then $\rep(D(G))$ is prime and therefore  by \cite[Proposition 4.5.]{gnn}  there is no triple (K,H,B), where K and H are normal subgroups of G that centralize each other, $(G, {e}) \neq (K,H) \neq ({e},G)$, $HK = G$, and B is a G-invariant bicharacter on $K \times H $such that the symmetric bicharacter $BB^{op}|_{(K\cap H)\times (K\cap H)}$ is nondegenerate.
\end{exmp}
\newpage

\bc
There is an isomorphism of associative algebras
\beqn
\ck(\cc) \ot k \simeq End_{\cc}(I(1)).
\eeqn
\ec

\bpf
Both algebras have the same block matrices but there is no canonical isomorphism.
\epf

\bc
Let $\cd$ be a modular category. Then $\cd=\cd_{ad}$ if and only if $\cd$ has no invertible elements. \blue{It follows also easily from the fact that $\cd_{ad}'=\cd_{pt}$.}
\ec
See [BK, Theorem 3.1.12] for the modular case. Moreover
\beqn
f_{\hat E_{A}}=\frac{\dim \cd}{(\dim A)^{2}}
\eeqn

Thus in the kernel of $C_{ad}$ are exactly those irreducible representations corresponding to the invertible objects of $\cc$.

Corollary for nilpotent modular categories if possible.

If $\cd$ is unimodular one has that:
\beqn
f_{\widehat E_{A}}=\frac{\dim \cd}{(\dim A)^{2}}.
\eeqn
It follows that $f_{{\widehat E_{A}}}=\dim \cd$ iff 
\section{The fusion subcategory generated by $I(1)$. A group homomorphism $U(\cz(\cc))\ra U(\cc)$.}
\subsection{On faithful gradings}
This generalize the results from \cite{bd}. Suppose that $\cc$ is a {\it faithfully} graded fusion category by the 
 group $G=U_{\cc}$. 
Therefore
$
\cc=\oplus_{g \in G}\cc_{g}.
$
Let $r_{g}$ be the regular element of the full abelian subcategory $\cc_{g}$. Thus
\beq
r_{g}=\sum_{X \in \irr(\cc_{g})}\fp(X)[X]
\eeq
In particular $r_{1}$ is the regular element of $\cc_{ad}$.
\md
Note that the regular element of $\cc$  is $r_{\cc}=\sum_{g\in G}r_{g}$. Suppose that $V \in \co(\cc_{h})$. Since inside $K_{0}(\cc)$ one has $[V]r=\fp(V)r$ it follows that 
\beq\label{r1}
[V]r_{g}=\fp(V)r_{gh}.
\eeq In particular $[V]r_{1}=\fp(V)r_{h}$.
\subsection{On the equivalence relations of a dominant tensor functor}
Let $F:\cc \ra \cd$ be a dominant tensor functor. In \cite{bubr} the authors defined an equivalence relation

Let $F : \cc \to \cd$ be a tensor functor between fusion categories, and denote by $I$ its right adjoint.

For $X \in \Lambda_\cc$ set
$$X^F = \{ [Y] \in \Lambda_\D \mid Y \,\mbox{is a factor of }\, F(X)\}$$ and
for $Y \in \Lambda_\D$, set $$Y_F = \{ [X] \in \Lambda_\C \mid Y \,\mbox{is a factor of }\, F(X)\}.$$

We have $X \in Y_F \iff Y \in X^F \iff X$ is a factor of $R(Y)$.

\subsection{On the forgetful functor} Apply the previous results for the forgetful functor $F:\cz(\cc)\ra \cc$. It is well known that the right adjoint is given by 

\beq\label{i1}
I(X)=\bigoplus_{i=1}^{r}V_{i}\ot X \ot V_{i^{*}}
\eeq
\blue{with the half braiding given as in Equation () of \cite{kir, eno, turviz}.}

By \cite[Proposition 4.3]{bubr} an equivalence class of the above equivalence relation $\approx^{F}$  consists of an entire coset of the fusion subcategory $<I(1)^{n}>$ relative to $\cz(\cc)$. 

Note that the forgetful image of the entire coset $I(1)^{n}$ consists of $\cc_{ad}$ the trivial component of the universal grading. It follows that 
There is a grading $\cz(\cc)$ defined by $M\in \cz(\cc)_{g}$ if $F(M)\in \co(\cc_{g})$.

\bt
The fusion subcategory of $\cz(\cc)$ generated by $I(1)$ coincides to the trivial component of the induced grading.
\et

\bpf
One needs to show that it does not split inside $\cc_{1}$. For two simple objects of $\cc$, by \cite[Lemma 4.7]{bubr}one has that $M\approx_{F}N$ if and only if $M$ is a factor of $N\ot F(I(1))^{n}$ for some $n \geq 0$. On the other hand Equation \eqref{I1} shows that $F(I(1))=C_{ad}$. Thus the fusion subcategory generate by $F(I(1))$ coincides to $\cc_{ad}$. Then $M\approx_{F}N$ if and only if they are in the same homogenous component.
\epf
From the previous theorem we get a canonical surjective group homomorphism $U(\cz(\cc))\ra U(\cc)$.
Let C be a fusion category and let $X \in  \irr(C)$. Recall that by the universal
property of U(X) there is a group homomorphism $ \phi_{X}:U_{X}\ra U_{\cc}$. It is defined by $\phi_{X}(t)=g$ if and only if $\cc(X)_{t}\subseteq\cc_{g}$.

\newpage
\section{}
Note that the first item of the previous theorem improves the statement of \cite[Corollary 1.3.7]{nfaith} which states that $K(A)$ is the unique maximal Hopf subalgebra contained in the intersection of the centers of all irreducible representations.

\beqn
\bigcap_{\pi \in \irr(A)} NZ(\pi)=K(A).
\eeqn
\blue{There is a compatibility on the tensor product of the Yetter-Drinfeld modules of the centers.}
\blue{see if the inequality is also true for noncommutative}
\section{Results from sonia}
\bt
The order of $g$ divides the order of $X$ if $X \in \co(\cc_{g})$.

If $\cc$ has $r$ homogenous generators then $U(\cc)$ also has $r$ generators, the degrees of the objects.
\et

\blue{cite this and extended it for braided fusion categories}.
\subsection{Normal fusion subcategories}
We call a  character of $K_{0}(\cc)$ \blue{positive} if it is a nonzero character which is a linear combination with nonnegative coefficients of the irreducible characters of $K_{0}(\cc)$.
\begin{lm}
Let $\cd$ be a normal fusion subcategory of a fusion category $\cc$. Then $\cc$ is the kernel of a positive character of $\cc$. 
\end{lm}
\bpf
Suppose that $\cd$ fits to the following short exact sequence of fusion categories:
\beq
\cd \hookrightarrow \cc \xra{F} \ce.
\eeq
One has an induced map $F_{!}:K_{0}(\cc) \ra K_{0}(\ce)$. For any irreducible character $\eta\in \ovr{K_{0}(\cc)}$ define $\hat{\eta}:=\eta \circ F_{!}$ the corresponding character of $K_{0}(\cc)$. It follows that 
\beq
\cd=\bigcap_{\eta\in \irr(K_{0}(\cc))} \ker_{\cc}(\hat{\eta}).
\eeq
Indeed, $X \in \cap_{\eta\in \irr(K_{0}(\cc))} \ker_{\cc}(\hat{\eta})$ if and only if $F(X)\in \cap_{\eta\in \irr(K_{0}(\cc))} \ker_{\ce}(\eta)$, i.e $X$ is a multiple of the identity of $\ce$.
\epf
\br
Use Hopf algebras to show that the converse is not true. There is a kernel that is not normal.
\er
\mdn
\blue{Is that true for any commutative fusion category?}
\mdn
\blue{See in \cite{DGNO} which relations are preserved in braided fusion categories}
\mdn
\red{Suppose that $H$ is a factorizable Hopf algebra. Then $LKer_H(M)$ is...}
\newpage
\red{Idea: Consider another generalization of $LKER_H(M)$ when the character ring is not commutative so one gets the right answer.}

\subsection{Relation with Hopf kernels}
Note that $\hat{d}$ is a trace on $K(A)$ Therefore $\hat{d}=\sum_{i \in \ce(d)}m_{d, i}\mu_i$ for some subset $\ce(d)\subset \{0, 1, \cdots , s\}$ and some nonnegative integers $m_{d, i}$.

There is a well defined linear function $\hat{ev}:C(A^*)\ra \widehat{C(A)}$ given by $d \mapsto \hat{d}$.
\bt
$d\in \ker_H(\ch)$ if and only if $\mtc{E}(d)\subset \ker_{\rep(H)}(\ch)$
\et
\red{If $C(H)$ is commutative then $\widehat{x_{ii}}$ is a morphism of $R(H)$; somehow unique: $\{coalgebras\} \ra X(C)$}
\mdn
\red{Is the union the entire kernel}?
\mdn
\red{Is the kernel a normal fusion subcategory?}

\blue{slightly nondegenerate how are they characterized; definition of modular categories.}\md
\blue{Let X be an object of C and let C[X] be the fusion subcategory of C generated
by X. The universal grading group of C[X] will be denoted by U(X).}
\section{Remarks, directions and questions}
\bne
\item generalization to the noncommutative case
\item generalizzato to the non semisimple case.
\item relation between the order of the kernel and that of the fusion subcategory generated.
\item As an example the category generated by the induced trivial module\; or  more general for fusion categories 
\item try  to use the formula for the centralizer from the mathematische zechfrift
\item singledg  might give an example of a modular tensor category that is prime but not single generated.
\item Find a counterexample for noncommutative rings! it does not make to much sense since $\ch(1)$ is not defined!
\item
\blue{The radical coincides with the algebraic radical, $\al^n=\ch(1)$, this is the center of a character; deduce more properties.} 
\item
\blue{normal tensor functors in terms of morphisms and kernels}
\item \blue{Frobenius property if $C(H)$ commutative}
\item The converse of the fact $U(\cc)$ cyclic implies that $\cc$ is generated by a single object, with $\cc$ nilpotent.
\item {Braided nilpotent what is the decomposition of the rings?}
\item \blue{What are the linear characters of the Grothendieck ring of $\Rep(A)$?}
\item {How can I characterize $\ck_0(<X>)$, $\ck_0(\cc(X))$ inside $\ck(\cc)$?}
\item {The interpretation of the class equation; and what a normal tensor functor means.}
\item Interpretation of the $Z(G/\ker \ch)\cong Z_{\ch}/\ker \ch$ in the fusion case.
\item Modular case; is this  a fusion subcategory.
\item If $\cc$ has a self dual object then $\fp(\cc)$ is even in the algebraic interring.
\item Replace $H^*$ by something to get that $ord(X)$ divides $multi(X)\fp(\cc)$.

\item A relation between the cardinal of the kernel and the dimension of the category generated by a simple object.
\item If the involution $b \mapsto \tilde{b}$ preserves $B_0$? if $\tilde{1}=1$, yes.
\item maybe it can be generalized to table Hopf algebras.
\item shorter paper only with kernels
\item ask Dimitry about quasi-semisimple Hopf
\item for the nonsemisimple case use example $GL_{p}(F_{p})$.
\item spherical fusion categories if needed.
\item pivotal structures in bijection with grouplike elements that inner the square of the antipode.
\item
Use the connection with $End(\id_{\cc})$ from Davydov's paper on jmp.
\item Generalize the divisibility theorems of Ito and that one relative to the center.
\item Folosind lemma de permutare Brauer cu $G=G(H)$ un subgroup al celui de grouplike elements. Numarul orbitelor actiunii pe $\irr(K)$ coincide cu nr orbitelor YD prin conjugare.
\item use species from sarah's thesis and generalize it 
\item \blue{Does it coincide to $\ker_{\cc}(X\ot X^{*})$? All the roots of unity form a subgroup of $\mathbb C^{*}$. Interpretation of the center in this case!}
\item \blue{
try to show that $RZ(\pi)$ as defined by Sonia coincides to $Z_{\ker \pi}$ as defined in this paper}
\item why do I include index and orders if not used?
\ene
      \newpage
\section{More properties of Grothendieck rings}
\bt
The order of $X$ divides $\fp(X)\mtr{multi}(X)$.
\et

\bc
Self dual objects implies $\fp(\cc)$ even.
\ec
\subsection{Coset decomposition for fusion categories}
\red{what are the principal eigenvalues; they coincide with those defined by the generating fusion category.}
\subsection{Double cosets in based rings}Let $\cd$ and $\ce$ be fusion subcategories of $\cc$. The space of double cosets $r_{\cd}\ck_0(\cc)r_{\ce}$ is spanned by $<E_j\;|j \in \ca_{\cd}\cap \ca_{\ce}>$.
\sub{Definition of the permutations}
\blue{ There is a permutation on the set of the indices $p$ such that $E_{i^*}=E_{p(i)}$.}
\mdn 
There is another permutation $I\mapsto i^*$ on the set of the indices such the $[X_i]^*=[X_{i^*}]$.

\subsection{Multiplicity between idempotents}
\bp
Let $\cc$ be a fusion category with a commutative Groethendieck ring $\ck_0(\cc)$. One has the following
\bne
\item
\beqn
[X]=\sum_{i=1}^{r_{\cc}}\frac{1}{a_i}m_{\cc}(E_i, [X^*])E_i
\eeqn
\item
\beq\lb{mid}
m_{\cc}(E_i, E_j)=a_j\delta_{i, \;p(j)}
\eeq
\ene
\ep
\bpf
For any object $X$ one has that
\bn{eqnarray*}
[X]&=&\sum_{i=1}^{r_{\cc}}m_{\cc}([X_i], [X])[X_i]=\sum_{i=1}^{r_{\cc}}m_{\cc}([X_i]^*, [X]^*)[X_i]\\&=&\sum_{i=1}^{r_{\cc}}\frac{1}{a_i}m_{\cc}(E_i, [X^*])E_i
\end{eqnarray*}

Then put $X=E_j$ and one obtains that
\beq\lb{e1}
m_{\cc}(E_i, E_j^*)=a_j\delta_{i,j}.
\eeq

Put also $X=E_j^*$ it follows that
\beq\lb{e2}
E_j^*=\sum_{i=1}^{r_{\cc}}\frac{1}{a_i}m_{\cc}(E_i, [E_j])E_i
\eeq

There is a permutation such that $X_{i}^*=X_i$ and a permutation such that $E_j^*=E_{p(j)}$.

Relation \ref{e1} shows that
\beq\lb{p1}
m_{\cc}(E_i, E_{p(j)})=a_j\delta_{i,j}.
\eeq
\blue{while the other shows the same thing.}
\epf
\br
For a Hopf algebra one has that $\ch_i^*=\ch_i\circ S$ and therefore $E_i^*=S_{H^*}(E_i)$. Therefore $p$ coincides with the permutation $^*$ given by $X_i^*=X_{i^*}$.
\er
\blue{Prove the same thing in general.}
\subsection{On the linear operators introduced on $\ck(\cc)$}
\subsubsection{The operator $c$}

Define an antilinear operator
\beq\lb{c}
c(\sum_{i}a_i[X_i])=\sum_{i}\overline{a_i}[X_i]
\eeq

One has that $c(\al x+\beta y)=\al c(x)+\overline{\beta}c(y)$. Also note that $c(xy)=c(x)c(y)$ for all $x, y \in \ck(\cc)$ and $c^2=\id$.

Indeed one has that
\beqn
xy=\sum_{k}(\sum_{i,j}N^k(i,j)x_iy_j)[X_k]
\eeqn
and therefore
\beqn
c(x)c(y)=\sum_{k}(\sum_{i,j}N^k(i,j)\overline{x_i}\overline{y_j})[X_k]=c(xy)
\eeqn

The operator $c$ induces a permutation $c$ on the set of the indices such that $c(E_i)=E_{c(i)}$.
\subsubsection{The operator $^o$}
Define also
\beq\lb{c}
(\sum_{i}a_i[X_i])^0=\sum_{i}\overline{a_i}[X_{i^*}]
\eeq
Note that $x^0=c(x)^*$. It follows that $(xy)^0=c(xy)^*=(c(x)c(y))^*=c(y)^*c(x)^*=y^0x^0$. Moreover ${x^0}^0=x$ and ${x^0}^*={x^*}^0=c(x)$.
\blue{Note that this the operator used by Witherspoon in \cite{reprint} and denoted also by $^*$.}
\subsubsection{The inner product}
As in \cite{repring} define
\beq\lb{ip}
<x, y>:=\sum_{i}x_i\overline{y_i}[X_i]
\eeq
Note that 
\beq
<x, y>=m_{\cc}(x, c(y))
\eeq
\blue{where $m_{\cc}$ is extended bilinearly on $\ck(\cc)\times \ck(\cc)$.}
Therefore 
\bn{eqnarray*}\lb{m1}
<x, yz>&=&m_{\cc}(x, c(yz))=m_{\cc}(x, c(z)c(y))\\&=&m_{\cc}(c(z)^*, c(y)x^*)=m_{\cc}(z^0, c(y)x^*)=m_{\cc}(z^0, c(x^0y))\\&=&<z^0,\;x^0y>
\end{eqnarray*}

Note also that
\beq
<x, y>=\overline{<y,x>}=<y^0,\;x^0>
\eeq

Same argument as in \cite{preprint} shows that $E_i^0=E_i$ and therefore 
\beqn
<E_i, E_j>=<1, E_iE_j>=\delta_{i,j}<1, E_j>
\eeqn
On the other hand
\beqn
<E_i, E_j>=m_{\cc}( E_i, c(E_j))=m_{\cc}( E_i, E_{c(j)})=\delta_{i,p(c(j))}a_i
\eeqn
which shows that $p(c(j))=\id$. This shows that $p=c^{-1}=c$, i.e $c(E_j)=E_j^*$ and therefore $\overline{\nu_{jl}}=\nu_{jl^*}$.
\mdn
\blue{I had some orthogonality relations in the noncommutative case}
\subsection{The conjugate value}
As in Proposition 3.1 of \cite{reprint} it follows that 
\beq\lb{conj}
\mu(x^0)=\overline{\mu(x)}
\eeq
for any object $x \in \ck(\cc)$. Therefore if $X=X_i$ is a simple object it follows that
\beq
\mu_i([X_{i^*}])=\overline{\mu([X_i])}
\eeq
\subsection{The formula for $N^h_{ij}$. Is this Verlinde's formula; generalized by Cohen for Hopf algebras.}
\bp
\blue{If $C(H)$ is commutative then $H$ is of Frobenius type}
\ep
\bpf
One has that 
\beq
e_{\ch}=\frac{1}{n}\sum_{(\Lam)}\ch(\Lam_1)S(\Lam_2)
\eeq
On the other hand 
\beq
\Lam=\sum_{d \in \Irr(A^*)}\eps(d)d
\eeq
where $d=\sum_{i}x^d_{ii}$.
Since $\ch(x_{ij})=\delta_{i,j}\ch(x_{ii})$ it follows that 
\beq
e_{\ch}=\frac{1}{n}\sum_{(d)}\eps(d)(\sum_{i}\ch(x_{ii})S(x_{ii}))
\eeq
Applying the algebra homomorphism $\omega_{\ch}:Z(H)\ra k$ given by $x \mapsto \frac{\ch(x)}{\ch(1)}$ it follows that
\beq
\frac{n}{n_{\ch}}=\sum_{(d)}\eps(d)(\frac{1}{\ch(1)}\sum_{i}\ch(x_{ii})\ch(S(x_{ii})))
\eeq
It follows that $\widehat{x_{ii}}=\mu_j$ for some $j$.
\epf
\section{On the induced ring morphism of a tensor functor}
Suppose that $\ck(\cc)$ is a commutative ring and $F:\cc \ra \cd $ is a tensor functor with a right adjoint $R: \cd\ra \cc$. Let $f=F_*$ and $r=R_*$. Then there is partition $\ca=\{\ca_s\}_s$ on the set $\{0, 1, \cdots , \mtr{rank}(\cd)-1\}$ such that
\beqn
F(E_i)=\sum_{j\in \ca_s}F_s
\eeqn
\br
It follows that $r(F_j)=E_{s(j)}$ where $s(j)$ is the unique index with $j \in \ca_s$.
\er
\bpf
Indeed by adjunction one has that
\beqn
m_{\cc}(r(F_j), E_t)=m_{\cd}(F_j, F(E_t))=\delta_{s(j), \;t^*}d_j
\eeqn

then the statement follows from equation \ref{mid}.
\epf
\begin{lm}Let $F:\cc \ra \cd$ be a tensor functor admitting a right adjoint $R:\cd\ra \cc$ and $f:=F_*:\ck_0(\cc)\ra \ck_0(\cd)$. One has that the morphism $f_*:\widehat{\ck_0(\cd)}\ra \widehat{\ck_0(\cc)}$ satisfies
\beqn 
f^*(\ker_{\cd}(F(X))\cong \ker_{\cc}(X)
\eeqn
for any object $X\in \co(\cc)$.
\end{lm}
\ncm{\rk}{\mtr{rank}}
\bp
A tensor functor $F:\cc \ra \cd$ is a normal tensor functor if and only if 
\beqn
\sum_{j\in \ca_s}d_j=\fp(R(1))a_s
\eeqn
for any $0 \leq s \leq \rk(\cc)-1$.
\ep

\bibliographystyle{ams}
\bibliography{mac-bob}
\ed